\documentclass{amsart}[10pt]

\usepackage{amscd,amssymb,latexsym,url,verbatim,graphicx,color,bm}
\usepackage{tikz,tikz-cd}
\usetikzlibrary{decorations.pathmorphing}

\usepackage{cases,amsmath}

\usepackage{txfonts}

\usepackage{tikz,tikz-cd}

\usepackage[dvips]{epsfig}

\title
{The first Cheeger constant of a simplex}

\author{Dmitry N. Kozlov}

\address{Department of Mathematics, University of Bremen, 28334
  Bremen, Federal Republic of Germany}

\email{dfk@math.uni-bremen.de}

\newtheorem{theorem}{Theorem}[section]
\newtheorem{df}[theorem]{Definition}
\newtheorem{thm}[theorem]{Theorem} 

\newtheorem{prop}[theorem]{Proposition}
\newtheorem{lm}[theorem]{Lemma}
\newtheorem{crl}[theorem]{Corollary}

\newtheorem{conj}[theorem]{Conjecture} 
\newtheorem{op}[theorem]{Open problem} 

\newcommand{\nin}{\noindent}
\newcommand{\pr}{\nin{\bf Proof.} }

\newcommand{\bo}{\partial}

\newcommand{\cmm}{{\mathcal CM}}

\newcommand{\da}{\Delta}

\newcommand{\sm}{\setminus}

\newcommand{\ti}{\tilde}

\newcommand{\dz}{{\mathbb Z}_2}
\newcommand{\sy} {\text{\rm sys}}
\newcommand{\csy}{\text{\rm csy}}

\newcommand{\bl}{{\bm{\lambda}}}
\newcommand{\bx}{\text{\rm box}}
\newcommand{\cor}{\text{\rm cor}}
\newcommand{\depth}{\text{\rm depth}}
\newcommand{\deff}{\text{\rm def}}
\newcommand{\val}{\text{\rm val}}

\numberwithin{equation}{section}
\numberwithin{figure}{section}
\numberwithin{table}{section}

\begin{document}

\begin{abstract}
The coboundary expansion generalizes the classical graph expansion to
the case of the general simplicial complexes, and allows the
definition of the higher-dimensional Cheeger constants $h_k(X)$ for an
arbitrary simplicial complex $X$, and any $k\geq 0$. In this paper we
investigate the value of $h_1(\da^{[n]})$ - the first Cheeger constant
of a simplex with $n$ vertices. It is known, due to the pioneering
work of Meshulam and Wallach, \cite{MW}, that 
\[\lceil n/3\rceil\geq h_1(\da^{[n]})\geq n/3, \textrm{ for all } n,\] 
and that the equality $h_1(\da^{[n]})=n/3$ is achieved when $n$ is
divisible by $3$.

\nin
Here we expand on these results. First, we show that
\[h_1(\da^{[n]})=n/3, \textrm{ whenever }n\textrm{ is not a power of }2.\] 
So the sharp
equality holds on a set whose density goes to $1$. Second, we show
that 
\[h_1(\da^{[n]})=n/3+O(1/n),\textrm{ when }n\textrm{ is a power of }2.\]  
In other words, as $n$ goes to infinity, the value
$h_1(\da^{[n]})-n/3$ is either $0$ or goes to $0$ very rapidly.

Our methods include recasting the original question in purely
graph-theoretic language, followed by a detailed investigation of
a~specific graph family, the so-called {\it staircase graphs}. These
are defined by associating a graph to every partition, and appear to
be especially suited to gain information about the first Cheeger
constant of a~simplex.
\end{abstract}

\maketitle

\section{Introduction}

The graph expanders are classical and well-studied mathematical
objects with many applications, see, e.g., the surveys \cite{HLW,
  Lu12}.  More recently, there have been different definitions of
higher-dimensional expanders, see~\cite{Lu14}. This paper is concerned
with the so-called {\it coboundary expanders}, which first made their
appearance in the paper by Linial and Meshulam, \cite{LiM}, and which
were later independently defined by Gromov, see~\cite{Gr}. Until now,
the major objective of the research on coboundary expansion has been
to find asymptotically good expanders, see, e.g., \cite{DK, LuM,LMM},
with computing the precise values of Cheeger constants playing the
secondary role. In this paper we deviate from this approach.

More specifically, the work we present here has a~twofold
purpose. Primarily, we are focused on taking the first step in the
general program of precise computation, or, at the very least, finding
sharp bounds for the higher Cheeger constants of standard simplicial
complexes. Currently, we do not even know the precise value of the
Cheeger constants for a~simplex.  In this paper, we attempt to change
that at least for the first Cheeger constant. To do that, we
reformulate the original questions for expansion in purely
graph-theoretical terms. Furthermore, in order to get the actual
estimates, we need to perform an in-depth analysis of certain graph
families.

Our second, more general purpose is to describe and to emphasize the
deep connection between the question of estimating the higher Cheeger
constants and questions in extremal graph and hypergraph theory. We
hope that this way the questions about coboundary expanders may gain
popularity and thus further progress on their understanding can be
achieved. In the conclusion of the paper we formulate several explicit
purely combinatorial conjectures.

Let us start by summarizing what is known about the Cheeger constants
of a~simplex with $n\geq 3$ vertices. First, a word about our 
notations. Usually writing $\da^n$ is reserved for the simplex 
of dimension $n$, that is the one having $n+1$ vertices. On the 
other hand, for an arbitrary set $V$ one uses the notation $\da^V$ 
to denote the simplex whose set of vertices is~$V$.  Since we also 
have the set notation $[n]:=\{1,\dots,n\}$, we find it consistent 
to use $\da^{[n]}$ to denote the simplex with $n$ vertices. 
The Cheeger constants $h_k(\da^{[n]})$ are then defined for all 
$0\leq k\leq n-2$, and so we are facing the task of determining
the numbers $h_0(\da^{[n]}),\dots,h_{n-2}(\da^{[n]})$. 

The $0$-th Cheeger constant is just the classical case and it is very
easy to calculate that $h_0(\da^{[n]})=\lceil(n+1)/2\rceil$, for all~$n$. 
On the other extreme, trivially one can see that
$h_{n-2}(\da^{[n]})=1$, for all~$n$.  Furthermore, it is not difficult 
to show, see Proposition~\ref{prop:pen}, that $h_{n-3}(\da^{[n]})=2$, 
for all~$n$. In general, we know, due to the work of Meshulam and
Wallach, see~\cite{MW}, that
\begin{equation}\label{eq:mw}
\lceil n/k\rceil\geq h_{k-2}(\da^{[n]})\geq n/k,
\end{equation}
for all $3\leq k\leq n$. Meshulam and Wallach also showed that the
lower bound is achieved when $k$ divides $n$.  On the other hand, we see
that the upper bound is sharp when $k=n-1$.

In this paper we are primarily concerned with the first Cheeger
constant $h_1(\da^{[n]})$. In this case $k=3$, and~\eqref{eq:mw}
specializes to $\lceil n/3\rceil\geq h_1(\da^{[n]})\geq n/3$, for all
$n\geq 3$, with equality $h_1(\da^{[n]})=n/3$ attained, whenever $n$
is divisible by~$3$.  Enhancing that information, we actually show
that $h_1(\da^{[n]})=n/3$, for all $n$, with a~definite exception of
the cases $n=4$ and $n=8$, and a~probable exception of the case when
$n$ is equal to other powers of~$2$. Furthermore, even when $n$ is a
power of $2$ we show that not only is $h_1(\da^{[n]})$ contained in
the interval between $n/3$ and $\lceil n/3\rceil$, but it actually
converges to $n/3$ very rapidly.  More specifically, we show that
$h_1(\da^{[n]})=n/3+O(1/n)$.

We finish this introductory chapter by describing briefly the plan of
the paper. In Section~\ref{sect:2} we recall the definition of the
coboundary expansion and the Cheeger constants. We then show how the
calculation of the first Cheeger constant can equivalently we
formulated as a graph-theoretic question. Section~\ref{sect:3} is the
core of the paper.  Here a family of graphs, which we call the {\it
  staircase} graphs is introduced and studied, computing all the
information which is relevant for the coboundary expansion. In
Section~\ref{sect:4} we apply the results of the previous section, both
to make precise calculation of the first Cheeger constant in the case
$n$ is not a~power of~$2$, as well as to derive sharp bounds in the
case $n$ is a~power of~$2$. We introduce the concept of a Cheeger
graph and find several of them realized as staircase graphs.  Finally,
in Section~\ref{sect:5} we state several open questions in extreme
graph and hypergraph theory, which are motivated by the coboundary
expansion. Section~\ref{sect:6} is the Appendix containing loose ends,
including the proof that $h_{n-3}(\da^{[n]})=2$, and recasting the coboundary 
computation of Wallach and Meshulam in the graph-theoretical language.

\newpage

\section{Setting up the board} \label{sect:2}

\subsection{The terminology of coboundary expanders} $\,$


\nin Let $X$ be a finite simplicial complex. In this paper, we shall
consider the associated chain and cochain complexes with
$\dz$-coefficients only, so we will suppress $\dz$ from the notations,
and simply write $C_*(X)$ and $C^*(X)$. Let now $\sigma$ be an arbitrary 
chain of~$X$, say $\sigma=\sigma_1+\dots+\sigma_d$, where $\sigma_i$ are 
generators indexed by the simplices of $X$, for $1\leq i\leq d$, and
$\sigma_i\neq\sigma_j$, whenever $i\neq j$. We set $\|\sigma\|:=d$ and
call this the {\it norm} of~$\sigma$. Dually, assume we have a~cochain
$c\in C^*(X)$, such that $c=c_1+\dots+c_d$, where $c_i$'s are
generators indexed by the distinct simplices of~$X$; each $c_i$ is the
characteristic function of a~$k$-simplex~$\sigma_i$. Then, we set 
$\|c\|:=d$, which we also call the norm of~$c$.


\begin{df}
For an arbitrary $k$-chain $\sigma$ we consider
$\min_\tau\|\sigma+\partial_*\tau\|$, where the minimum is taken over all $(k+1)$-chains $\tau$. 
We call that number the {\bf systolic norm} of
$\sigma$, and denote it by $\|\sigma\|_\sy$. Furthermore, a~{\bf
  systolic form} of $\sigma$ is any
$\tilde\sigma=\sigma+\partial_*\tau$, such that
$\|\ti\sigma\|=\|\sigma\|_\sy$. We let $\sy(\sigma)$ denote the set of
all systolic forms of $\sigma$.  A chain is called a~{\bf systole} if
$\|\sigma\|=\|\sigma\|_\sy$.


Dually, assume $c$ is a~$k$-cochain.  We call
$\min_d\|c+\partial^* d\|$, where the minimum is taken over all $(k-1)$-cochains 
the {\bf cosystolic norm} of $c$, and
denote it by $\|c\|_\csy$.  Let $\csy(c)$ denote the set of all
cosystolic forms of $c$.  A~{\bf cosystolic form} of $c$ is any
$\tilde c=c+\partial^* d$, such that $\|\ti c\|=\|c\|_\csy$.  A chain
is called a~{\bf cosystole} if $\|c\|=\|c\|_\csy$.
\end{df}

The cosystolic norm of a cochain can be quite difficult to compute in general.

\begin{df}
Assume we are given a simplicial complex $X$. For any $k$-cochain $c$ of
$X$, which is not a~coboundary, the {\bf coboundary expansion} of $c$ is
\[\|c\|_{\exp}:=\|\partial^* c\|/\|c\|_\csy.\]  
The $k$-th {\bf Cheeger constant} of $X$ is then
\begin{equation}\label{eq:dfh}
h_k(X):=\min_{c\neq\partial^* d}\|c\|_{\exp}.
\end{equation}
\end{df}

Clearly, in~\eqref{eq:dfh} we might as well restrict ourselves to
cosystoles, when taking the minimum. Finally, when $c$ is a cosystole
such that $\|c\|_{\exp}=h_k(X)$, then we shall call $c$ a~{\it Cheeger
  cosystole}.

\subsection{Simplicial complex of cut-minimal graphs} $\,$

\nin Let us now introduce some graph terminology in order to give an
alternative definition of the first Cheeger constant. We shall use the
notation $G=(V,E)$, meaning that the graph $G$ has the set of vertices
$V$ and the set of edges~$E$. For any two, not necessarily disjoint,
subsets $A,B\subset V$, we set $E(A,B):=|\{(v,w)\in E\,|\,v\in
A,\,w\in B\}|$, and $NE(A,B):=|\{(v,w)\notin E\,|\,v\in A,\,w\in B\}|$,
so $E(A,B)+NE(A,B)=|A|\cdot |B|$.

\begin{df}
 A graph $G=([n],E)$ is called {\bf cut-minimal} if for any proper
 subset $S\subset[n]$ we have
\begin{equation}\label{eq:cmg}
|E(S,[n]\sm S)|\leq |NE(S,[n]\sm S)|.
\end{equation}
In other words, at most half of the $|S|(n-|S|)$ potential edges
connecting vertices from $S$ to vertices from $[n]\sm S$ belong
to~$G$.

We call the cut $(S,[n]\sm S)$ {\bf perfect} if equality is achieved
in~\eqref{eq:cmg}.
\end{df}

In particular, the valencies of vertices of a~cut-minimal graph with $n$ 
vertices can be at most $\lfloor(n-1)/2\rfloor$. All the graphs which are 
shown on Figure~\ref{fig:1} are cut-minimal. The way we think about the 
condition~\eqref{eq:cmg} is as follows. Imagine we are given a graph $G$ and 
we are allowed to split the vertex set $[n]$ into two parts: $S$ and $[n]\sm S$. 
We take all the potential edges between these two parts, and think of them as
a~{\it cut} $C$. We are now allowed change $G$ by inverting the {\it
  being the edge of $G$} relationship within $C$. In other words, we
obtain a new graph 
by keeping all the edges in $G$
which are outside of $C$, removing all the edges of $G$ which are in
$C$ and adding as edges all the non-edges of $G$ which are in $C$. The
graph is then cut-minimal if no such operation can decrease the number
of edges of $G$; which explains our choice of terminology.

Note that removing some edges from  a cut-minimal graph will certainly yield a cut-minimal
graph again. Following the general ideology of combinatorial topology,
see~\cite{book}, this observation leads to a definition of a natural
combinatorial simplicial complex. 

\begin{df}
Let us fix $n\geq 2$. The abstract simplicial complex $\cmm(n)$ is defined
as follows:
\begin{itemize}
\item the vertices are indexed by unordered pairs $\{i,j\}$,
  $i,j\in[n]$, $i\neq j$;
\item the set of vertices forms a simplex of $\cmm(n)$ if and only if
  the corresponding graph is cut-minimal.
\end{itemize}
\end{df}

We see that $\cmm(2)$ is empty, $\cmm(3)$ has $3$ vertices and no
edges, and the complex $\cmm(4)$ has $6$ vertices and $3$ disjoint
edges. The complex $\cmm(5)$ is more interesting. It has dimension $3$
and its $f$-vector is $(10,45,100,10)$. In particular, $\cmm(5)$ has
$10$ vertices and a~full $1$-skeleton. It can be obtained from a~full
$2$-skeleton by deleting $20$ triangles and adding $10$
tetrahedra. Its maximal simplices are these $10$ tetrahedra, together
with $60$ triangles. It can be shown by direct inspection, using
a~combination of techniques from \cite{book}, that $\cmm(5)$ is
homotopy equivalent to a~wedge of $54$ spheres of dimension~$2$.

In general, clearly $\cmm(n)$ has $\binom{n}{2}$ vertices, for all
$n\geq 3$.  Furthermore, it is non-pure for all $n\geq 5$. It would be
interesting to understand more the simplicial structure or topology of
these complexes.  For example, the dimension of $\cmm(n)$ is obtained
by subtracting $1$ from the maximal number of edges which a
cut-minimal graph may have. This number has been computed precisely in
the upcoming work~\cite{KR}.

\subsection{The graph-theoretic definition of the first Cheeger constant} $\,$

\nin Assume now we are given a simplicial complex $X$. Its
$1$-skeleton $G:=X^{(1)}$ is a graph, whose set of vertices is
$V:=X(0)$ and whose set of edges is $E:=X(1)$. Here we follow very
handy notations of Linial and Meshulam, \cite{LiM}, by letting $X(k)$
denote the set of all $k$-simplices of $X$.

The edges of this graph $G$ are in 1-to-1 correspondence with the
generators of the group of 1-cochains $C^1(X)$: associate to each edge
$e$ its characteristic cochain which evaluates to $1$ on $e$ and to
$0$ on all other edges. For simplicity we identify each edge with the
associated characteristic cochain.  Since we are working over~$\dz$,
the arbitrary cochains can be identified with the {\it sets of edges}
of $X$, or, which is the same, with the subgraphs of~$G$.

In the same way, the vertices of $G$ are in 1-to-1 correspondence with the
generating $0$-cochains of $X$, and sets of vertices of $G$ are in
1-to-1 correspondence with arbitrary $0$-cochains.  Taking the
coboundary has graph-theoretic translation too. Given an arbitrary
$0$-cochain $c$ corresponding to a set of vertices $S$, its coboundary
is the $1$-cochain which corresponds to the edge set $E(S,V\sm S)$.
The norm of the $0$-cochain is $|S|$, and the norm of the $1$-cochain
is $|E(S,V\sm S)|$.

\begin{prop}
The correspondence above restricts to a 1-to-1 correspondence between
the sets of cosystoles and cut-minimal graphs.
\end{prop} 
\pr Being a~cosystole means that addition of any coboundary will not
increase norm.  This is the same as to say that the corresponding
graph contains at most half of the edges in the induced cut.  \qed

\vspace{5pt}

\nin The following definition associates a certain number to an
arbitrary graph.

\begin{df}\label{df:tvwu}
Assume we are given a graph $G=([n],E)$. For each edge $e=(v,w)\in E$,
we set $t(e):=\sum_{u\in [n],u\neq v,w} \tau_e(u)$, where the numbers
$\tau_e(u)$ are defined as follows:
\begin{equation}\label{eq:tvwu}
\tau_e(u):=
\begin{cases}
1,& \text{ if }(v,u),(w,u)\notin E;  \\
1/3, & \text{ if } (v,u),(w,u)\in E;\\ 
0,& \text{ otherwise. } \\
\end{cases} 
\end{equation}
We now set $h(G):=\sum_{e\in E}t(e)/|E|$. 
\end{df}
\nin Definition~\ref{df:tvwu} can alternatively be phrased as follows.
Let $T(G)$ denote the set of all ``triangles'' which contain an odd
number of edges from $G$, i.e.,
\begin{equation}\label{eq:tg}
T(G):=\{(v,e)\,|\,v\in V,\,e=(w,u)\in E,\, v\notin e,\,
\left|\{(v,w),(v,u),(u,w)\}\cap E\right|\text{is odd}\}.
\end{equation}
We have $\sum_{e\in E}t(e)=|T(G)|$. This is because,
by~\eqref{eq:tvwu}, if a~triangle from $T(G)$ has one edge from $G$,
then this edge gives a contribution $1$ to the sum $\sum_{e\in
  E}t(e)$, and if a~triangle from $T(G)$ has three edges from $G$,
then each of these edges gives a contribution $1/3$ to that sum.
We therefore have the alternative formula
\begin{equation}\label{eq:tg2}
|E|\cdot h(G)=|T(G)|.
\end{equation}

\nin We are now ready to give a graph-theoretic description of the
first Cheeger constant of a~simplex.

\begin{prop}
For any $n\geq 3$ we have
\[h_1(\da^{[n]})=\min_G h(G),\]
where the minimum is taken over all cut-minimal graphs $G$ with $n$
vertices.
\end{prop}
\pr By definition, the constant $h_1(\da^{[n]})$ is equal to
$\min_c\|\bo^* c\|/\|c\|$, where the minimum is taken over all
cosystoles~$c$. As mentioned above, being a~cosystole precisely
corresponds to cut-minimal graphs, and computing the value $h(G)$ is
exactly the same as computing $\min_c\|\bo^* c\|/\|c\|$. \qed

\vspace{5pt}

\nin It is rather straightforward to extend this description to the
first Cheeger constant of an arbitrary simplicial complex.


\section{Staircase graphs} \label{sect:3}

\subsection{Terminology of partitions} $\,$

\nin
A {\it partition} $\bl$ is any ordered tuple of positive integers
$(\lambda_1,\dots,\lambda_t)$, such that
$\lambda_1\geq\lambda_2\geq\dots\geq\lambda_t$. In such a case, we
always set the default values $\lambda_q:=0$, for all $q>t$.
The {\it Ferrers diagram} of a partition $\bl=(\lambda_1,\dots,\lambda_t)$ 
is the arrangement of square boxes in $t$ rows, such that the boxes are left-justified, 
the first row is of length $\lambda_1$, the second row is of length $\lambda_2$, and so on.
When referring to the individual boxes in the diagram, we shall count both rows and
columns starting with $1$, counting rows from top to bottom and counting columns from
left to right. 

To abbreviate our writing, we shall use
the power notation for the diagram, i.e., using formal powers to denote multiple 
parts of the same cardinality, for example:
$(3^{(2)},2^{(3)},1)=(3,3,2,2,2,1)$.
For $\bl=(\lambda_1,\dots,\lambda_t)$, we set
$|\bl|:=\sum_{k=1}^t\lambda_k$.  We also set $\bx(\bl):=\lambda_1+t$,
say $\bx(3^{(2)},2^{(3)},1)=9$.

Given a partition $\bl=(\lambda_1,\dots,\lambda_t)$, a {\it conjugate
  partition} $\bl^*=(\mu_1,\dots,\mu_m)$ is defined as follows: we set
$m:=\lambda_1$, and for every $1\leq k\leq m$, we set $\mu_k$ to be
equal to the maximal index $i$ such that $\lambda_i\geq k$. In
particular of course $\mu_1=t$. In terms of the Ferrers diagram we just
switch rows and columns of $\bl$. For an arbitrary partition $\bl$, we 
have $(\bl^*)^*=\bl$, $|\bl^*|=|\bl|$, and $\bx(\bl)=\bx(\bl^*)$.
As an example, we have $(3^{(2)},2^{(3)},1)^*=(6,5,2)$.

For an arbitrary $t\geq 1$, we let $\cor(t)$ denote the partition
$(t,t-1,\dots,2,1)$. We clearly have $\bx(\cor(t))=2t$, $|\cor(t)|=t(t+1)/2$, and
$\cor(t)^*=\cor(t)$.

\begin{df}\label{df:depth}
Let $\bl=(\lambda_1,\dots,\lambda_t)$ be an arbitrary partition. The
{\bf depth} of~$\bl$, denoted $\depth(\bl)$, is the maximal number $d$
such that the Ferrers diagram of $\cor(d)$ is contained in the Ferrers
diagram of~$\bl$.
\end{df}

Alternatively, the depth of $\bl$ can be described as the unique value
$d$, such that
\begin{enumerate}
\item[(1)] $\lambda_1\geq d,\lambda_2\geq d-1,\dots,\lambda_k\geq d-k+1,
\dots,\lambda_d\geq 1$,
\item[(2)] there exists $1\leq k\leq d+1$, such that $\lambda_k=d-k+1$.
\end{enumerate}

A convenient way to think about $\depth(\bl)$ is to notice that it is equal to the minimal
number of rows and columns which will cover the entire Ferrers diagram of $\bl$,
or, expressed algebraically, we have
\[\depth(\bl)=\min_{0\leq k\leq t}(k+\lambda_{k+1}),\]
where we use the convention $\lambda_{t+1}=0$.

Of course, we have $\depth(\cor(t))=t$.

\subsection{The definition of staircase graphs} $\,$

\nin
The following family of graphs is central to our approach.

\begin{df}
Assume we are given a partition $\bl$ and an~integer $n$, such that
$n\geq\bx(\bl)$.  The {\bf staircase graph}
$G_n(\lambda_1,\dots,\lambda_t)=G_n(\bl)$ is defined as follows:
\begin{itemize}
\item the set of vertices of $G_n(\bl)$ is a~disjoint union $V\cup
  W\cup U$, where $V=\{v_1,\dots,v_l\}$, with $l=\lambda_1$,
  $W=\{w_1,\dots,w_t\}$, and $U=\{u_1,\dots,u_r\}$, with $r=n-l-t$;
\item for each $1\leq i\leq l$, $1\leq j\leq t$, the vertices $v_i$
  and $w_j$ are connected by an~edge if $\lambda_j\geq i$; all other
  pairs of vertices are not connected by an edge.
\end{itemize}
\end{df}

\begin{figure}[hbt]

  \input{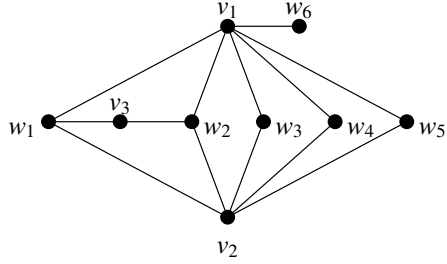}  

\caption{The graph $G_{9}(3,3,2,2,2,1)$.}
\label{fig:sg}
\end{figure}

\begin{figure}[hbt]

  \input{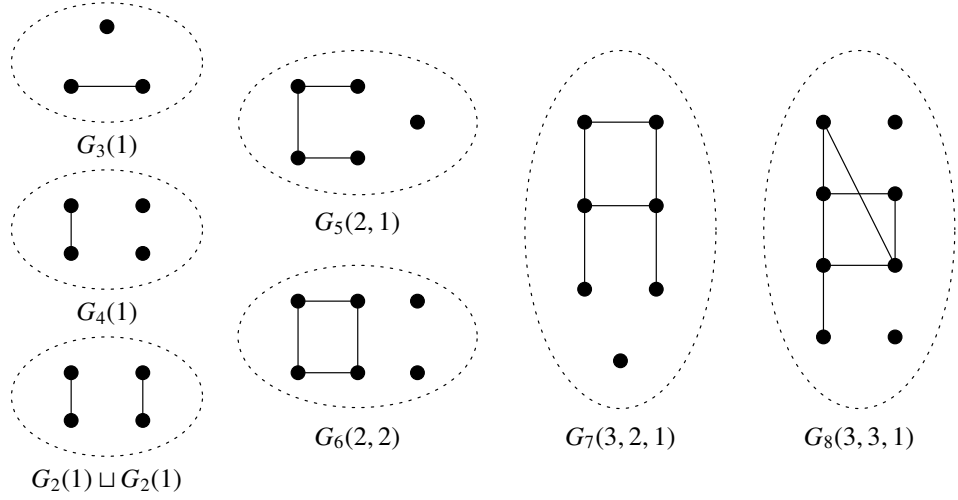}  

\caption{Up to isomorphism, these are all Cheeger graphs on $n$
  vertices, for $3\leq n\leq 8$.}
\label{fig:1}
\end{figure}

In particular, we see that $G_n(\bl)$ is always bipartite, $V$ and $W$
can be taken as two sides of the bipartition, and vertices of $U$ are
isolated. Note that $l,t\neq 0$, whereas $r$ might be $0$; this will
happen if $n=\bx(\bl)$.
Clearly, we also have $|V(G_n(\bl))|=n$, and $|E(G_n(\bl))|=|\bl|$.
Figure~\ref{fig:sg} shows a~staircase graph, several further examples 
can be found on Figure~\ref{fig:1}.
Finally, note that $G_n(\bl)$ is isomorphic to $G_n(\bl^*)$.

\subsection{Structure theory of staircase graphs} $\,$

\nin
Assume we are given a partition $\bl=(\lambda_1,\dots,\lambda_t)$, and
$1\leq k\leq t$, and $1\leq m\leq\lambda_1$. 

\begin{df} 
For arbitrary index sets $I\subseteq[t]$ and $J\subseteq[\lambda_1]$, such
that $|I|=k$, and $|J|=m$, let $B_{\bl}(I,J)$ denote the total number
of boxes in the Ferrers diagram of $\bl$, which are either contained in
one of the $k$ rows indexed by $I$, or contained in one of the $m$
columns indexed by~$J$, but not in both.
\end{df}

If we consider a cut of $G_n(\bl)$ with the vertices corresponding to the 
row and column indices from $I\cup J$, then $B_\bl(I,J)$ is precisely the number of edges
across that cut.

\begin{df} \label{df:legal}
Assume we are given a partition $\bl$, and an integer $n\geq\bx(\bl)$.
Assume $\bl^*=(\mu_1,\dots,\mu_q)$. The partition $\bl$ is called {\bf
  legal with respect to} $n$ if the following three conditions are
satisfied:
\begin{enumerate}
\item[(1)] for all $1\leq k\leq t$, we have $B_\bl([k],\emptyset)=
\sum_{i=1}^k\lambda_i\leq k(n-k)/2$;
\item[(2)] for all $1\leq m\leq q$, we have $B_\bl(\emptyset,[m])=
\sum_{j=1}^m\mu_j\leq m(n-m)/2$;
\item[(3)] $|\bl|\leq\depth(\bl)(n-\depth(\bl))/2$ .
\end{enumerate}
\end{df}

Let us say a few words on the intuition behind Definition~\ref{df:legal}.
Conditions (1) and (2) make sure that the cut-minimality holds when we cut
off the vertices corresponding to the first $k$ rows or the first $m$
columns. Condition~(3) is rather concerned with the cuts where we choose $k$
first rows and $\lambda_{k+1}$ first columns. This condition could actually be 
strengthened to require that 
\[|\bl|-k\lambda_{k+1}\leq(k+\lambda_{k+1})(n-k-\lambda_{k+1})/2,\]
for all $k=1,\dots,t-1$. We do not need this strengthening here and find it 
technically simpler to work with the condition in Definition~\ref{df:legal}.

It turns out that legality of a partition has the following strong implication.

\begin{lm}
Assume $n$ is a natural number, and a~partition $\bl$ is legal with
respect to~$n$.  Choose arbitrary index sets $I\subseteq[t]$ and
$J\subseteq[q]$, such that $|I|+|J|\leq n/2$, and set $k:=|I|$,
$m:=|J|$. We have
\begin{equation}\label{eq:bij}
B_\bl(I,J)\leq(k+m)(n-k-m)/2.
\end{equation}
\end{lm}
\pr Without loss of generality we can shift all the rows upwards and
all the columns to the left. If after this they cover the entire
Ferrers diagram of $\bl$, then $B_\bl(I,J)\leq|\bl|$ and
$k+m\geq\depth(\bl)$.  We then get
\[B_\bl(I,J)\leq\depth(\bl)(n-\depth(\bl))/2\leq(k+m)(n-k-m)/2,\]
where the first inequality is given by condition (3) of
Definition~\ref{df:legal}, and the second inequality follows from the
fact that $\depth(\bl)\leq k+m\leq n/2$.

If, on the other hand, the Ferrers diagram is not covered completely,
we have
\[
\begin{split} 
B_\bl(I,J)&=(\lambda_1+\dots+\lambda_k)+(\mu_1+\dots+\mu_m)-2km\\
&\leq k(n-k)/2+m(n-m)/2-2km \\
&=(k+m)(n-k-m)/2-km\leq (k+m)(n-k-m)/2,
\end{split}
\]
where the first inequality follows from conditions (1) and (2) of
Definition~\ref{df:legal}.  
\qed

\vspace{5pt}

Clearly, if a~partition $\bl$ is legal with respect to some $n$, and
$n'\geq n$, then $\bl$ is also legal with respect to~$n'$. This observation 
motivates the following definition.

\begin{df} 
For an arbitrary partition $\bl$, we let $N(\bl)$ denote the minimal 
natural number with respect to which $\bl$ is legal.
\end{df}

For example, one can compute that $N_{3,3,1}=8$, and $N_{6,5,2}=13$. 
Of course, we have $N_\bl=N_{\bl^*}$.

Given $\bl=(\lambda_1,\dots,\lambda_t)$, we set 
\[N_r(\bl):=\max_{1\leq k\leq t}
\left(k+\left\lceil\frac{2(\lambda_1+\dots+\lambda_k)}{k}\right\rceil\right),\]
\[N_d(\bl):=\depth(\bl)+\left\lceil\frac{2|\bl|}{\depth(\bl)}\right\rceil.\]
The following lemma gives us a precise formula for computing $N(\bl)$.

\begin{lm}
For an arbitrary partition $\bl$, we have
\[N(\bl)=\max\left\{N_r(\bl),N_r(\bl^*),N_d(\bl)\right\}.\]
\end{lm}
\pr Simply rewrite the inequalities of Definition~\ref{df:legal}.
\qed

\vspace{5pt}

\nin
As an example, for any $t\geq 1$, we get
\[N_d(\cor(t))=t+\left\lceil
2\cdot\frac{t(t+1)}{2}\cdot\frac{1}{t}\right\rceil=t+(t+1)=2t+1\]
\[N_r(\cor(t))=\max_{1\leq k\leq t}
\left(k+\left\lceil\frac{2(t+\dots+(t-k+1))}{k}\right\rceil\right)=
\max_{1\leq k\leq t}(k+(2t-k+1))=2t+1,\]
so, since $\cor(t)=\cor(t)^*$, we conclude that $N(\cor(t))=2t+1$.

\begin{lm}
If partition $\bl$ is legal with respect to some number $n$, then the
graph $G_n(\bl)$ is cut minimal.
\end{lm}
\pr Cutting the set $[n]$ into the subsets $S$ and $[n]\sm S$, such
that $|S|\leq n/2$, is the same as choosing subsets $I$ and $J$, with
$k=|I|$ and $m=|J|$, such that $n-\bx(\bl)\geq |S|-k-m\geq 0$. Then
$B_\bl(I,J)$ is the number of edges across the cut, and we have
$(k+m)(n-k-m)/2\leq|S|(n-|S|)/2$, so~\eqref{eq:bij} implies the
cut-minimality.  \qed

\begin{df}
For any $\bl=(\lambda_1,\dots,\lambda_k)$,
$\bl^*=(\mu_1,\dots,\mu_m)$, we set
\[|\bl^2|:=\dfrac{1}{2}\left(\sum_{i=1}^k\lambda_i^2+\sum_{j=1}^m\mu_j^2\right),\]
and furthermore, we set
\[h(\bl):=N(\bl)-\frac{2|\bl^2|}{|\bl|}.\]
\end{df}

We remark that
our notation $|\bl^2|$ is the special case of
$|\bl^p|:=\left(\sum_{i=1}^k\lambda_i^p+\sum_{j=1}^m\mu_j^p\right)/2$,
which for $p=1$ also gives our notion~$|\bl|$.
For future reference, for an arbitrary partition $\bl$, we set 
\[\deff(\bl):=h(\bl)-N(\bl)/3=\frac{2}{3}\left(N(\bl)-
\frac{3|\bl^2|}{|\bl|}\right),\]
which we call the {\it deficiency} of $\bl$.

Taking $\bl=\cor(t)$ as a~specific example, we can see 
\[|\cor(t)^2|=1^2+\dots+t^2=t(t+1)(2t+1)/6,\]
which implies
\[h(\cor(t))=2t+1-2\cdot\frac{t(t+1)(2t+1)}{6}\cdot\frac{2}{t(t+1)}=
2t+1-\frac{2}{3}(2t+1)=\frac{2t+1}{3},\] 
and hence $\deff(\cor(t))=0$.

\begin{lm}\label{lm:hbl}
We have $h(\bl)=h(G_{N(\bl)}(\bl))$.
\end{lm}
\pr In our notations, we have 
\[t(e_{ij})=(N(\bl)-\lambda_1-\mu_1)+(\mu_1-\mu_j)+(\lambda_1-\lambda_i)=
N(\bl)-\lambda_i-\mu_j,\]
where we recall Definition~\ref{df:tvwu}. Therefore
\begin{equation}
\label{eq:hgnbl}
h(G_{N(\bl)}(\bl))=\dfrac{1}{|\bl|}\sum_{i,j}t(e_{i,j})=
N(\bl)-\dfrac{1}{|\bl|}\sum_{i,j}(\lambda_i+\mu_j)=h(\bl),
\end{equation}
where the sum is taken over all $i$ and $j$, which correspond to boxes
in the Ferrers diagram of $\bl$. Note, that the last equality follows from the fact 
that $\sum_{i,j}\lambda_i=\sum_i\lambda_i^2$ (this is because for each $i$ the number 
of summands $\lambda_i$ on the left hand side is equal to the number of boxes in 
the $i$th row of the Ferrers diagram of $\bl$, and this number is of course precisely 
$\lambda_i$), and analogously $\sum_{i,j}\mu_j=\sum_j\mu_j^2$.
\qed

\subsection{The partition $c\bl$.} \label{ssect:clambda} $\,$

\nin The staircase partitions can be blown up using the following
simple operation.

\nin
\begin{df}
Given a partition $\bl=(\lambda_1,\dots,\lambda_t)$, and a natural
number $c$, we set
\[c\bl:=(c\lambda_1,\dots,c\lambda_1,\dots,c\lambda_t,\dots,c\lambda_t)=
(\ti\lambda_1,\dots,\ti\lambda_{ct}),\] 
where ${\ti\lambda}_q=c\lambda_{\lceil q/c\rceil}$, for all $1\leq q\leq ct$.
\end{df}

\nin
The next lemma relates the data associated to the partition $c\bl$ to
the data associated to the partition $\bl$.

\begin{lm}\label{lm:many}
For an arbitrary partition $\bl$ and an arbitrary natural number $c$
we have the following equalities:
\begin{equation}\label{eq:cbl-depth}
\depth(c\bl)=c\cdot\depth(\bl),
\end{equation}
\begin{equation}\label{eq:cbl-val}
|c\bl|=c^2|\bl|,
\end{equation}
\begin{equation}\label{eq:cbl-cont}
|(c\bl)^2|=c^3\cdot|\bl^2|
\end{equation}
\begin{equation}\label{eq:cbl-nd}
N_d(c\bl)=c\cdot N_d(\bl)
\end{equation}
and inequalities
\begin{equation}\label{eq:cbl-nr}
N_r(c\bl)\leq c\cdot N_r(\bl),
\end{equation}
\begin{equation}\label{eq:cbl-n}
N(c\bl)\leq c\cdot N(\bl)
\end{equation}
\begin{equation}\label{eq:cbl-h}
h(c\bl)\leq c\cdot h(\bl)
\end{equation}
\end{lm}
\pr We start by showing~\eqref{eq:cbl-depth}. For brevity, set
$d:=\depth(\bl)$, and let us show that $cd$ satisfies conditions (1)
and~(2) in the Definition~\ref{df:depth} for
$c\bl=(\ti\lambda_1,\dots,\ti\lambda_{ct})$. First, for all
$1\leq\ti\lambda_q\leq cd$, we see that
\[\ti\lambda_q=c\lambda_{\lceil q/c\rceil}\geq c(d-\lceil q/c\rceil+1)\geq
cd-(q+c-1)+c=cd-q+1,\] where we used the inequality $c\cdot\lceil
q/c\rceil\leq q+c-1$.  This verifies condition (1). Second, assume
$\lambda_k=d-k+1$, for some $1\leq k\leq d+1$. Then, $1\leq ck-c+1\leq
dc+1$, and we have
\[\ti\lambda_{ck-c+1}=c\lambda_k=cd-ck+c=cd-(ck-c+1)+1,\]
which verifies condition~(2).

The equalities \eqref{eq:cbl-val}, \eqref{eq:cbl-cont}, and
\eqref{eq:cbl-nd}, are direct computations which we leave to the
reader.

Let us show the inequality~\eqref{eq:cbl-nr}. We have
\[N_r(c\bl)=\max_{1\leq k\leq ct}\left(k+\left\lceil
\frac{2(\ti\lambda_1+\dots+\ti\lambda_k)}{k}\right\rceil\right).\] 
To start with, for $k=cm$, for $1\leq m\leq t$, we have
\begin{multline*}
cm+\left\lceil\frac{2(\ti\lambda_1+\dots+\ti\lambda_{cm})}{cm}\right\rceil=
cm+\left\lceil\frac{2(c^2\lambda_1+\dots+c^2\lambda_m)}{cm}\right\rceil=
cm+\left\lceil\frac{2c(\lambda_1+\dots+\lambda_m)}{m}\right\rceil\leq \\ \leq
cm+c\cdot\left\lceil\frac{2(\lambda_1+\dots+\lambda_m)}{m}\right\rceil\leq
cN_r(\bl),
\end{multline*} 
where the penultimate inequality follows from the fact that $\lceil
cx\rceil\leq c\lceil x\rceil$, whenever $c$ is an integer.

Next, consider the special case $1\leq k\leq c-1$. Then, we have
\[k+\left\lceil\frac{2(\ti\lambda_1+\dots+\ti\lambda_{k})}{k}\right\rceil=
k+\left\lceil\frac{2kc\lambda_1}{k}\right\rceil=k+2c\lambda_1<c(1+2\lambda_1)\leq
c\cdot N(\bl).
\]

Assume now that $k=cp+r$, where $1\leq p\leq t-1$, $1\leq r\leq
c-1$. Set $\ti r:=r/c$, so $0<\ti r<1$.  We have
\begin{multline*}
cp+r+\left\lceil\frac{2(\ti\lambda_1+\dots+\ti\lambda_{cp+r})}{cp+r}\right\rceil=
cp+r+\left\lceil\frac{2(c^2(\lambda_1+\dots+\lambda_p)+cr\lambda_{p+1})}{cp+r}
\right\rceil=\\ c\left(p+\ti
r+\frac{1}{c}\left\lceil\frac{2c(\lambda_1+\dots+\lambda_p+\ti
  r\lambda_{p+1})} {p+{\ti r}}\right\rceil\right)\leq c\left(p+\ti
r+\left\lceil\frac{2(\lambda_1+\dots+\lambda_p+\ti r\lambda_{p+1})}
{p+{\ti r}}\right\rceil\right)
\end{multline*}
Set $s:=\lambda_1+\dots+\lambda_p$, and note that
$\lambda_{p+1}<s/p$. We now claim that
\begin{equation}\label{eq:conv}
p+\ti r+\left\lceil\frac{2(s+\ti r\lambda_{p+1})}{p+\ti r}\right\rceil\leq 
(1-\ti r)\left(p+\left\lceil\frac{2s}{p}\right\rceil\right)
+\ti r\left(p+1+\left\lceil\frac{2(s+\lambda_{p+1})}{p+1}\right\rceil\right).
\end{equation}
Clearly this would finish our proof of~\eqref{eq:cbl-nr}, since the
right hand side of \eqref{eq:conv} is bound above by $(1-\ti r)
N_r(\bl)+\ti r N_r(\bl)=N_r(\bl)$.  On the other hand, the inequality
\eqref{eq:conv} is equivalent to
\begin{equation}\label{eq:srl1}
\left\lceil\frac{2(s+\ti r\lambda_{p+1})}{p+\ti r}\right\rceil\leq
(1-\ti r)\left\lceil\frac{2s}{p}\right\rceil +\ti r
\left\lceil\frac{2(s+\lambda_{p+1})}{p+1}\right\rceil.
\end{equation}
A direct calculation, using the fact that $\lambda_{p+1}<s/p$ shows that
\begin{equation}\label{eq:srl2}
\frac{2(s+\ti r\lambda_{p+1})}{p+\ti r}\leq 
(1-\ti r)\frac{2s}{p}+\ti r\frac{2(s+\lambda_{p+1})}{p+1}.
\end{equation}
We now apply $\lceil -\rceil$ to both sides of~\eqref{eq:srl2}, and
use the fact that $\lceil x+y\rceil\leq\lceil x\rceil+\lceil y\rceil$
to verify~\eqref{eq:srl1}.

Finally, the inequalities \eqref{eq:cbl-n} and \eqref{eq:cbl-h}
are both obtained by a direct substitution.
\qed


\begin{crl}\label{crl:eqn}
Assume that we are given a~partition $\bl$, such that
$N(\bl)=N_d(\bl)$, then we actually have the equalities
$N(c\bl)=c\cdot N(\bl)$, and $h(c\bl)=c\cdot h(\bl)$.
\end{crl}
\pr
We have 
\[c\cdot N(\bl)\geq N(c\bl)\geq N_d(c\bl)=c\cdot N_d(\bl)=c\cdot N(\bl),\]
where the first inequality is~\eqref{eq:cbl-n}, the second inequality
is the definition of $N(-)$, the penultimate equality
is~\eqref{eq:cbl-nd}, and the last equality is the assumption of the
corollary.  This shows that $N(c\bl)=c\cdot N(\bl)$, and the equality
$h(c\bl)=c\cdot h(\bl)$ is an immediate consequence.  \qed


\section{Applications of staircase graphs}\label{sect:4}

\subsection{Exact value of the first Cheeger constant for the simplex whose
number of vertices is not a~power of $2$} $\,$

\vspace{5pt}

\nin Let us set $h(n):=h_1(\da^{[n]})=\min_G h(G)$, for $n\geq 3$,
where the minimum is taken over all cut-minimal graphs $G$ with $n$
vertices.

\begin{thm} \label{thm:lmw}
(Meshulam-Wallach bound, \cite{MW})

\nin For any $n$ we have
\[\lceil n/3\rceil\geq h(n)\geq n/3.\]
In particular, if $3$ divides $n$ then we have $h(n)=n/3$.
\end{thm}

We provide a~short write-up of the proof of the lower bounds of
Theorem~\ref{thm:lmw} using our notations in
subsection~\ref{ssect:6.4}.  Extending this result, our next theorem
shows that the lower bound of Theorem~\ref{thm:lmw} is true for the
vast majority of the values of $n$. We will give two proofs of the
following theorem, one here using a~direct computation, and one in the
appendix as Corollary~\ref{crl:hkm}.

\begin{thm}\label{thm:main1}
Assume $n=c(2t+1)$, such that $t\geq 1$, then $h(n)=n/3$.
\end{thm}
\pr Consider the partition $\bl:=c\cdot\cor(t)$. Lemma~\ref{lm:many}
implies that $|\bl|=c^2 t (t+1)/2$, $\depth(\bl)=ct$, and
$N_d(\bl)=c(2t+1)=n$. Since $N_d(\cor(t))=2t+1=N(\cor(t))$,
Corollary~\ref{crl:eqn} implies that $N(\bl)=n$ and $h(\bl)=n/3$.  

Let us now consider the corresponding staircase graph $G=G_n(\bl)$.
By Lemma~\ref{lm:hbl}, we have $h(G)=h(\bl)=n/3$. It follows
immediately from Theorem~\ref{thm:lmw} that $h(G)=h(n)=n/3$.  \qed

\begin{crl}
If $n\neq 2^\alpha$, then we have $h(n)=n/3$.
\end{crl}
\pr If $n$ is not a power of $2$, then it can be written 
as $n=c(2t+1)$, where $t\geq 1$.
\qed

\vskip5pt

\nin
The following terminology seems natural, allowing us to talk about the
graphs that are optimal with respect to the first Cheeger constant.

\begin{df}
A graph $G$ with $n$ vertices is called a {\bf Cheeger graph} if
$h(G)=h(n)$.
\end{df}

\nin We have already found a~family of Cheeger graphs.

\begin{crl}
Assume $c$ and $t$ are arbitrary natural numbers. Set $\bl:=c\cdot\cor(t)$,
and recall that $N(\bl)=c(2t+1)$. Then the graph $G_{c(2t+1)}(\bl)$ is
a~Cheeger graph with $c(2t+1)$ vertices.
\end{crl}
\pr This follows directly from the proof of Theorem~\ref{thm:main1}.
\qed

\vspace{5pt}

\nin Note, how the special case $t=1$ yields graphs considered in the
previous work of Meshulam and Wallach. In general, it would be interesting
to describe the set of all Cheeger graphs.


\subsection{Bounds for the first Cheeger constant for the simplex whose
number of vertices is a~power of $2$} $\,$

\nin Since $n\geq 3$, the first relevant power of $2$ is $2^2=4$, in
which case we have $h(4)=2$. Furthermore, specific examples, and in
case $n=8$, exhaustive case analysis show that
\[\arraycolsep=1.4pt
\begin{array}{rcccl}
h(8) -\dfrac{8}{3}  &   = &\dfrac{4}{21}   &\approx &0.19,\\[0.4cm]
h(16)-\dfrac{16}{3} &\leq &\dfrac{8}{93}   &\approx &0.086,\\[0.4cm]
h(32)-\dfrac{32}{3} &\leq &\dfrac{16}{381} &\approx &0.042.
\end{array}
\]

\nin In general we have the following upper bound.

\begin{thm} \label{thm:main2}
Assume $n=2^d$, for some $d\geq 3$, then we have
\begin{equation}\label{eq:2db}
h(n)-n/3\leq\frac{4n}{3(n^2-8)}=\frac{4}{3n}+\frac{32}{3n(n^2-8)}=\frac{1}{3}
(4n^{-1}+32n^{-3}+256n^{-5}+\dots).
\end{equation}
\end{thm}
\pr Set $t:=n/4\geq 2$. Set
\[\bl:=((2t-1)^{(2)},(2t-3)^{(2)},\dots,3^{(2)},1).\]
 Note that the conjugate partition is given by
\[\bl^*=(2t-1,(2t-2)^{(2)},\dots,2^{(2)}).\] 
We have
\[|\bl|=|\bl^*|=2\sum_{k=1}^t(2k-1)-1=2t^2-1,\]
and
\[|\bl^2|=2\sum_{k=1}^{2t-1}k^2+(2t-1)^2-1=\frac{1}{3}t(8t^2-5).\]

\nin A direct check shows that
$N_d(\bl)=N_r(\bl)=N_r(\bl^*)=N(\bl)=4t=n$. Hence the staircase graph
$G_{N(\bl)}(\bl)=G_n({\bl})$ is cut-minimal.

Substituting the obtained values into the formula for $h(\bl)$, we
obtain the following calculation:
\[h(\bl)=4t-\frac{2t(8t^2-5)}{3(2t^2-1)}=\frac{2t(4t^2-1)}{3(2t^2-1)}.\]
We conclude that
\[\deff(\bl)=h(\bl)-n/3=\frac{2t}{3(2t^2-1)}=\frac{4n}{3(n^2-8)}=
\frac{4}{3n}+\frac{32}{3n(n^2-8)},\]
showing the inequality~\eqref{eq:2db}.
\qed

\vspace{5pt}

\nin The next corollary is immediate.

\begin{crl}
We have $\lim_{n\to\infty}(h(n)-n/3)=0$.
\end{crl}


\section{Conjectures and open problems}\label{sect:5}

\nin We know that the following Conjecture~\ref{conj:A} is true for
$\alpha=1$, $2$, and $3$.

\begin{conj} \label{conj:A}
We have $h(n)>n/3$, for all $n=2^\alpha$, $\alpha\geq 1$.
\end{conj}

\nin In all of our examples, the constant $h(G)$ for optimal graphs
$G$ never had a~contribution coming from a triangle with all $3$ edges
in $G$, in other words, the second line of~\eqref{eq:tvwu} was
invoked. We conjecture that this holds in general.

\begin{conj} \label{conj:B}
All Cheeger graphs are triangle-free.
\end{conj}

\nin
We actually believe that a stronger statement is true.

\begin{conj} \label{conj:B2}
All Cheeger graphs are bipartite.
\end{conj}

\nin The next conjecture is very daring, and would clearly imply
Conjectures~\ref{conj:B} and~\ref{conj:B2}.

\begin{conj} \label{conj:C}
All Cheeger graphs except for $G_2(1)\sqcup G_2(1)$ can be
represented as staircase graphs.
\end{conj}

\nin We finish this section with two open problems, which are probably
rather hard, but which might help to stimulate further research.

\begin{op}
Classify all Cheeger graphs, for $n\geq 9$.
\end{op}

\begin{op}
Determine the topology of the simplicial complexes $\cmm(n)$, for
$n\geq 6$.
\end{op}

\nin We mention, that recently, see~\cite{Me16}, the asymptotics of
these simplicial complexes, and, more generally, of the simplicial
complexes of $k$-cosystoles, has been understood.



\newpage

\section{Appendix}\label{sect:6}

\nin This humble section contains some facts and elementary proofs
which we feel would be useful to fix in writing for future reference.




\subsection{Blowing up the graphs} $\,$

\nin
Let us now generalize the blowing of partitions, which we did in
subsection~\ref{ssect:clambda}, to blowing up arbitrary graphs.

\begin{df}
Let $G=(V,E)$ be an arbitrary graph, and let $c$ be any natural
number.  We let $cG=(\widetilde V,\widetilde E)$ to be the graph
defined as follows:
\begin{itemize}
\item we set $\widetilde V:=V\times[c]$;
\item for any $v,w\in V$, and $i,j\in[c]$, we have $((v,i),(w,j))\in
  \widetilde E$ if and only if $(v,w)\in E$.
\end{itemize}
\end{df}

\nin In particular, we have $|V(cG)|=c|V(G)|$ and
$|E(cG)|=c^2|E(G)|$. To connect this to our partition notations, we
note that for any partition $\bl$ we have
$cG_n(\lambda)=G_{cn}(c\lambda)$.

\begin{prop} \label{prop:hcg}
For an arbitrary graph $G$ and any natural number $c$, we have
\[h(cG)=c\cdot h(G).\]
\end{prop}
\pr It follows directly from the definition in~\eqref{eq:tg}, that
$T(cG)=c^3\cdot T(G)$. Hence~\eqref{eq:tg2} implies that
\[h(cG)=\frac{T(cG)}{E(cG)}=\frac{c^3\cdot|T(G)|}{c^2\cdot|E(G)|}=c\cdot h(G).\]
\qed

\vspace{5pt}

\nin The next theorem is the main result of this subsection.

\begin{thm}\label{thm:ccm}
Assume $G=(V,E)$ is a~cut-minimal graph, and $c$ is an arbitrary
natural number. Then, the graph $cG$ is also cut-minimal.
\end{thm}
\pr Before proceeding with a formal argument, we would like to give a
informal idea of how the proof goes. If $cG$ is not cut-minimal then
it must have a ``bad'' cut. This cut cannot nicely go around the blown
up vertices, as in these cuts the number of edges simply changes
proportionally and the original graph was cut-minimal. So the bad cut
must cut at least one of the blown up vertices. Now, shifting vertices
between the two parts of the cut within the blown up vertices changes
the number of edges which cross the cut linearly, while the total
number of potential edges in the cut changes along a concave
function. This means that one of these changes will yield a bad cut
again, and so eventually we will get a bad cut which does not cut any
of the blown up vertices, leading to a~contradiction.

Let us now make this argument rigorous.  For simplicity of notations,
we set $cV:=V(cG)$ and $cE:=E(cG)$.  Let us take an arbitrary proper
subset $S\subset cV$. Assume first that $S= T\times[c]$, for some
proper subset $T\subset V$.  Then
\[|E(S,cV\sm S)|=c^2\cdot|E(T,V\sm T)|\leq c^2\cdot|NE(T,V\sm T)|=
|NE(S,cV\sm S)|,\] where the sets of edges and non-edges are always
taken in the appropriate graphs.  This verifies the
condition~\eqref{eq:cmg} for the set $S$ and graph~$cG$.

Assume that $cG$ is not cut-minimal. It follows from the previous
paragraph that we can pick $S\subset cV$, such that the
condition~\eqref{eq:cmg} is not satisfied, and there does not exist
any proper subset $T\subset V$, such that $S=T\times[c]$. This means
that we can pick $v\in V$, such that both sets $A:=S\cap(\{v\}\times
[c])$ and $B:=(cV\sm S)\cap(\{v\}\times [c])$ are non-empty. Clearly
$\{v\}\times [c]$ is a~disjoint union of $A$ and $B$.

Set $S^-:=S\sm A$, and $S^+:=S\cup B$. We shall use a~concavity
argument to show that condition~\eqref{eq:cmg} is not satisfied for at
least one of the sets $S^-$ and $S^+$. For the short-hand notations we
set $s:=|S|$, $a:=|A|$, $b:=|B|$, $n:=|V|$, $s^+:=|S^+|=s+b$,
$s^-:=|S^-|=s-a$, $e:=E(S,cV\sm S)$, $e^+:=E(S^+,cV\sm S^+)$, and
$e^-:=E(S^-,cV\sm S^-)$. Note that $a+b=c$. Let $w$ be any vertex from
$\{v\}\times[c]$, and set $\beta:=|E(w,S)|$, and $\gamma:=|E(w,cV\sm
S)|$. Note that these numbers do not depend on the choice of $w$. Also
$\beta+\gamma=c\cdot\val(v)$, but we will not need that. Note that
moving such a vertex $w$ between $A$ and $B$ changes the number of
the edges of the graph which cross the cut by $\beta-\gamma$, so we
have
\begin{equation}\label{eq:ccm1}
e^-=e+a(\beta-\gamma)\,\,\text{ and }\,\,e^+=e+b(\gamma-\beta),
\end{equation}
which yields
\begin{equation}\label{eq:ccm2}
ae^++be^-=ce.
\end{equation}

Assume both $S^+$ and $S^-$ satisfy condition~\eqref{eq:cmg}. This
means that $e^+\leq s^+(cn-s^+)/2$ and $e^-\leq
s^-(cn-s^-)/2$. Combining these with~\eqref{eq:ccm2} we get
\begin{equation}\label{eq:ccm3}
\frac{a}{a+b}s^+(cn-s^+)+\frac{b}{a+b} s^-(cn-s^-)\geq e.
\end{equation}
The function $f(x)=x(cn-x)$ is concave, which means that
\[\frac{a}{a+b}f(s+b)+\frac{b}{a+b}f(s-a)\leq 
f\left(\frac{a}{a+b}(s+b)+\frac{b}{a+b}(s-a)\right)=f(s).\] 
This translates to
\[\frac{a}{a+b}s^+(cn-s^+)+\frac{b}{a+b} s^-(cn-s^-)\leq s(cn-s),\]
which together with~\eqref{eq:ccm3} contradicts to the fact that
condition~\eqref{eq:cmg} is not satisfied for~$S$.

Repeating this argument we can modify $S$ until it has a form
$T\times[c]$, while the condition~\eqref{eq:cmg} is still not
satisfied. This clearly contradicts the first paragraph of this proof,
so we are done.  \qed

\vspace{5pt}

\nin
We can now derive a~generalization of \eqref{eq:cbl-h} as a~simple
corollary of Theorem~\ref{thm:ccm}.

\begin{crl} \label{crl:hkm}
For any $c\geq 1$, and any $n\geq 3$, we have
\begin{equation}\label{eq:hcn}
h(cn)\leq c\cdot h(n).
\end{equation}
\end{crl}
\pr Take any Cheeger graph $G$ with $n$ vertices. We have
$h(G)=h(n)$. The graph $cG$ has $cn$ vertices, and, by
Theorem~\ref{thm:ccm} it is cut-minimal. It follows that $h(cG)\geq
h(cn)$. On the other hand, by Proposition~\ref{prop:hcg}, we have
$h(cG)=c\cdot h(G)=c\cdot h(n)$, hence~\eqref{eq:hcn} follows. \qed

\vspace{5pt}

\nin Note, that~\eqref{eq:hcn} implies that if $h(n)=n/3$, then
$h(cn)=cn/3$, for all natural numbers~$c$. This yields another, and
simple proof of Theorem~\ref{thm:main1}, since we can limit ourselves
to the analysis of the staircase graphs associated to $\cor(t)$,
which, in turn, is rather straightforward.


\subsection{Computing the penultimate Cheeger constant of a simplex}$\,$

\nin As promised, we now provide a~simple argument for precise
computation of $h_{n-3}(\da^{[n]})$. In this case, the upper
Meshulam-Wallach bound is realized.

\begin{prop}\label{prop:pen}
For any $n\geq 3$, we have $h_{n-3}(\da^{[n]})=2$.
\end{prop}
\pr By definition, we have $h_{n-3}(\da^{[n]})=\min_{c}\|\bo^*
c\|/\|c\|_\csy$, where the minimum is taken over all $c\in
C^{n-3}(\da^{[n]})$, $c\neq\bo^* f$. The group of cochains
$C^{n-3}(\da^{[n]})$ is generated by characteristic cochains of
simplices of codimension~$2$, i.e., simplices with $n-2$ vertices.
For all $1\leq k,l\leq n$, $k\neq l$, let $c_{kl}$ denote the
characteristic $(n-3)$-cochain of the $(n-3)$-simplex $[n]\sm\{k,l\}$;
that is, $c_{kl}$ evaluates to $1$ on that simplex and it evaluate to
$0$ on all other $(n-3)$-simplices. We have $c_{kl}=c_{lk}$. Each
cochain $c\in C^{n-3}(\da^{[n]})$ has a~unique presentation as a~sum
$c_{k_1 l_1}+\dots+c_{k_t l_t}$, such that
$\{k_i,l_i\}\neq\{k_j,l_j\}$, for all $i\neq j$. 

For $i=1,\dots,n$, let $d_i\in C^{n-2}(\da^{[n]})$ denote the
characteristic $(n-2)$-cochain of the $(n-2)$-simplex $[n]\sm\{i\}$.
We clearly have 
\begin{equation}\label{eq:dckl}
\bo^* c_{kl}=d_k+d_l.
\end{equation}

We have a bijection between the generators $c_{kl}$ and edges of
a~complete graph on $n$ vertices~$K_n$. If we extend this bijection to
the one between $d_i$'s and vertices of $K_n$, then the coboundary
equation~\eqref{eq:dckl} translates to taking the boundary of an edge
in that graph. Note, that~\eqref{eq:dckl} means that for all $c\in
C^{n-3}(\da^{[n]})$, not just the characteristic ones, we know that
$\|\bo^* c\|$ must be even, since all the cancellations happen in
pairs. The fact that $H^{n-3}(\da^{[n]})=0$, implies that if $c$ 
is not a~coboundary, then it is not a~cocycle, i.e., $\bo^* c\neq 0$,
so $\|\bo^* c\|\geq 2$.

Finally, let $f_{klm}\in C^{n-4}(\da^{[n]})$ denote the characteristic
cochain of the simplex $[n]\sm\{k,l,m\}$, for all $1\leq k,l,m\leq n$,
$k\neq l\neq m$. We clearly have
\begin{equation}\label{eq:dfklm}
\bo^* f_{klm}=c_{kl}+c_{km}+c_{lm}.
\end{equation}

Let us now pick an arbitrary non-zero cosystole $c\in
C^{n-3}(\da^{[n]})$, and write $c=c_{k_1 l_1}+\dots+c_{k_t l_t}$, with
$\{k_i,l_i\}\neq\{k_j,l_j\}$, for all $i\neq j$. Assume that not all
the numbers in the set $\{k_1,\dots,k_t,l_1,\dots,l_t\}$ are distinct.
Then, without loss of generality, we can assume that $k=k_1=k_2$. The
equation~\eqref{eq:dfklm} implies that $c_{kl_1}+c_{kl_2}+c_{l_1l_2}$
is a~coboundary. Adding this expression to $c$ would decrease the
norm, which contradicts the fact that we picked $c$ to be a~cosystole.
Thus, we can assume that all the numbers in the set
$\{k_1,\dots,k_t,l_1,\dots,l_t\}$ are distinct. Since
$\bo^*c=d_{a_1}+\dots+d_{a_t}+d_{b_1}+\dots+d_{b_t}$, we conclude that
$\|\bo^* c\|=2t$, and hence $\|c\|_{\exp}=2t/t=2$.

On the other hand, non-zero $(n-3)$-dimensional cosystoles clearly
exist. For example, $c_{12}$ is such a~cosystole. Indeed,
$\bo^*c_{12}=d_1+d_2\neq 0$, so $c_{12}$ is not a~coboundary, i.e.,
$\|c_{12}\|_\csy\geq 1$. On the other hand
$1=\|c_{12}\|\geq\|c_{12}\|_\csy$, hence $\|c_{12}\|_\csy=1$.

We conclude that in dimension $n-3$, all non-zero cosystoles are
in fact Cheeger cosystoles, and that $h_{n-3}(\da^{[n]})=2$.
\qed

\subsection{The proof of the lower bound in Meshulam-Wallach theorem 
using our notations}
\label{ssect:6.4} $\,$

\nin We restrict ourselves to proving the lower bound from
Theorem~\ref{thm:lmw}, since the upper bound is improved by other
results in this paper. The argument below follows closely the ideas of
the original coboundary computation by Meshulam and Wallach,
\cite{MW}. Still we find it instructive, and potentially useful, to
phrase it in our elementary language.



\vspace{5pt}

\nin {\bf Proof of the lower bound in Theorem~\ref{thm:lmw}.} Assume
we are given a cut-minimal graph $G=(V,E)$, such that $|V|=n$. Let $M$
denote the set of all ordered pairs $(v,e)$, where $v\in V$, and
$e=(w,u)\in E$, such that $v\notin e$, and the number of edges of $G$
among $(v,w)$, $(v,u)$, and $(w,u)$ is odd, in other words, set
\[M:=\{(v,e)\,|\,v\in V,\,e=(w,u)\in E,\, v\notin e,\,
\left|\{(v,w),(v,u),(u,w)\}\cap E\right|\text{ is odd}\}.\] Comparing
with \eqref{eq:tg}, we immediately see that
$|M|=3\cdot|T(G)|$, hence~\eqref{eq:tg2} implies that
\begin{equation}\label{eq:ap0}
|M|=3\cdot|E|\cdot h(G).
\end{equation}

\begin{figure}[hbt]

  \input{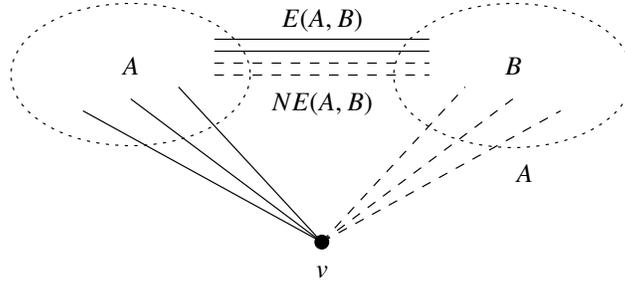}  

\caption{Proof of the lower Meshulam-Wallach bound.}
\label{fig:2}
\end{figure}

\nin On the other hand, for a~fixed $v\in V$, let $M_v$ denote the set
of edges $e$, such that $(v,e)\in M$; clearly $|M|=\sum_{v\in
  V(G)}|M_v|$. We now show that $|M_v|\geq |E|$, for all $v\in V$.
Set $A:=\{w\in V\,|\,(v,w)\in E\}$, and $B:=\{w\in V\,|\,(v,w)\notin
E\}$, see Figure~\ref{fig:2}. If $A=\emptyset$, then $M_v=E(B,B)=E$,
so we might as well assume $A\neq\emptyset$, and $(A,V\sm A)$ is
a~proper cut.  We have
\begin{equation}\label{eq:ap1}
M_v=E(A,A)\sqcup E(B,B)\sqcup NE(A,B),
\end{equation}
which is to be compared with
\begin{equation}\label{eq:ap2}
E=E(A,A)\sqcup E(V\sm A,V\sm A)\sqcup E(A,V\sm A)=
  E(A,A)\sqcup E(B,B)\sqcup E(A,V\sm A).
\end{equation}
By definition of $A$, we have $NE(A,V\sm A)=NE(A,B)$. Hence, the
cut-minimality condition~\eqref{eq:cmg}, when applied to the cut
$(A,V\sm A)$, yields $E(A,V\sm A)\leq NE(A,B)$. Together
with~\eqref{eq:ap1} and~\eqref{eq:ap2}, this proves $|M_v|\geq |E|$,
for all $v\in V$.  Since $|M|=\sum_{v\in V(G)}|M_v|$, we get $|M|\geq
n\cdot|E|$.  Together with~\eqref{eq:ap0} this yields $h(G)\geq n/3$,
for all cut-minimal graphs $G$, and hence $h(n)\geq n/3$.  \qed

Note, that one can also see from our proof of the lower bound in
Theorem~\ref{thm:lmw}, that the sharp bound $h(n)=n/3$ is achieved by
a~cut-minimal graph $G$ if and only if for every non-isolated vertex
$v$ the corresponding cut $(A,V\sm A)$ is perfect. This observation
gives us a~quick-and-dirty argument for the strict inequality
$h(8)>8/3$. Indeed, the size of $A$ is a valency of $v$, so if $v$ is
not isolated, it is equal to $1$, $2$, or $3$, as $G$ is
cut-minimal. If $(A,V\sm A)$ is a~perfect cut, then $|A|$ must be
even, otherwise $|A|\cdot(8-|A|)$ would have been odd. This means that
$|A|=2$, and all non-isolated vertices of $G$ have valency $2$. The
graph $G$ is a~disjoint union of isolated vertices and cycles, and
$h(G)\geq n-4$ for such graphs. Here this means $h(G)\geq 4>8/3$.

\section*{Acknowledgment}
\nin The author would like to thank Roy Meshulam for drawing his
attention to the exciting topic of coboundary expanders. He also would
like to thank Eva-Maria Feichtner for useful discussions. This
research has been supported by the grant No.\ 1261 of the
German-Israeli Foundation for Scientific Research and Development,
DFG-ANR Grant ``DISCMAT'', as well as by University of Bremen.


\end{document}